\documentclass[12pt]{extarticle}
\usepackage{amsmath, amsthm, amssymb, hyperref, color}
\usepackage[capitalize,noabbrev]{cleveref}
\usepackage{graphicx}
\usepackage{caption}
\usepackage{subcaption}
\usepackage{mathtools}
\usepackage{faktor}
\usepackage{enumerate}
\usepackage[all]{xypic}
\usepackage{verbatim}
\usepackage{chemfig, chemnum}
\usepackage{tikz}
\usepackage{tikz-cd}
\usepackage{caption}
\usepackage{enumitem}
\usepackage{subcaption}
\oddsidemargin \evensidemargin
  \hypersetup{
  	colorlinks=true,
 	linkcolor=blue,
 	filecolor=blue,      
 	citecolor=blue,
 	urlcolor=blue,
 }

\theoremstyle{theorem}
\newtheorem{theorem}{Theorem}
\newtheorem{proposition}[theorem]{Proposition}
\newtheorem{lemma}[theorem]{Lemma}
\newtheorem{corollary}[theorem]{Corollary}
\theoremstyle{definition}
\newtheorem{definition}[theorem]{Definition}

\newtheorem{remark}[theorem]{Remark}
\newtheorem{conjecture}[theorem]{Conjecture}

\newtheorem{example}[theorem]{Example}

\DeclareMathOperator{\val}{val}
\DeclareMathOperator{\PZ}{\mathrm{PZ}}
\DeclareMathOperator{\GL}{GL}

\DeclareMathOperator{\End}{\mathrm{End}}

\newcommand{\RR}{\mathbb{R}}
\newcommand{\QQ}{\mathbb{Q}}

\newcommand{\ZZ}{\mathbb{Z}}
\newcommand{\Ocal}{\mathcal{O}}

\newcommand{\Bcal}{\mathcal{B}}

\newcommand{\minplus}{\,\underline{\oplus}\,}
\newcommand{\maxplus}{\,\overline{\oplus}\,}
\newcommand{\mintimes}{\,\underline{\odot}\,}
\newcommand{\maxtimes}{\,\overline{\odot}\,}

\title{\textbf{Orders and Polytropes: \\ 
		Matrix Algebras from Valuations}}

\author{ 
	Yassine El Maazouz,
	Marvin Anas Hahn, Gabriele Nebe, \\
	Mima Stanojkovski,
	Bernd Sturmfels}

\date{}
\begin{document}
	
	\maketitle
	
	\begin{abstract}
		\noindent We apply tropical geometry to study matrix algebras over a field with valuation.
		Using the shapes of min-max convexity, known as polytropes, we revisit the graduated orders introduced by Plesken and Zassenhaus.
		These are classified by the polytrope region. We advance the ideal theory of graduated orders by introducing their
		ideal class polytropes. This article emphasizes examples and computations. It offers first steps in the geometric combinatorics of endomorphism rings of configurations in affine buildings.
	\end{abstract}	
	
	\section{Introduction}
	
	Let $K$ be a field with a surjective discrete valuation
	$\val: K \to \ZZ \cup \{\infty\} $. We fix 
	$p \in K$ satisfying~${\rm val}(p) = 1$.
	The \emph{valuation ring}  $\Ocal_K$ is the
	set of elements in $K$ with non-negative valuation. This is a local  ring
	with maximal ideal
	$\langle p \rangle = \{ x \in \Ocal_K : \val(x) > 0\}$.
	In our examples,
	$K = \QQ$ is the field of rational numbers, with the $p$-adic valuation
	for some prime~$p$.
	
	We write $K^{d \times d}$ for the ring of $d \times d$ matrices
	with entries in $K$. The map ${\rm val}$ is applied coordinatewise
	to matrices and vectors. For example, if $K = \QQ$ with $p=2$, then the vector
	$x = ( 8/7 , 5/12 , 17)$ has ${\rm val}(x) = (3 , -2 , 0)$.
	In what follows, we often take $ X = (x_{ij})$ 
	to be a $d \times d$ matrix with nonzero entries in $K$.
	In this case,  ${\rm val}(X) = ({\rm val}(x_{ij}))$ 
	is a matrix in $\ZZ^{d \times d}$.

	Fix any square matrix $M = (m_{ij})$ in $\ZZ^{d\times d}$.
	This paper revolves around the interplay between the following two objects
	associated with $M$,
	one algebraic and the other geometric:
	\begin{enumerate}
		\item the  set $\Lambda_M=\{ X \in K^{d \times d} : \val(X) \geq M \}$,
		an $\Ocal_K$-lattice in the vector space $K^{d \times d}$;
		\item the set $Q_M = \{ u\in \RR^d/\RR{\bf 1} : u_i - u_j \leq m_{ij} 
		\textup{ for }  1 \leq i,j \leq d \} $, where ${\bf 1} = 
		(1,\ldots,1)$.
	\end{enumerate}
	This interplay is strongest and most interesting when $\Lambda_M$ is closed under multiplication.
	In this case, $\Lambda_M$ is a non-commutative ring of matrices.
	Such a ring is called an \emph{order} in~$K^{d\times d}$. 
	The  quotient space  $ \RR^d/\RR{\bf 1} \simeq \RR^{d-1}$ is
	the usual setting for tropical geometry \cite{joswig, maclagan}.
	Note that $Q_M$ is a convex polytope in that space. It is also tropically convex,
	for both  the min-plus algebra and the max-plus algebra. Following
	\cite{joskul, tran}, we use the term {\em polytrope} for~$Q_M$.
	
	\begin{example}
		\label{ex:rhombic}
		For $d=4$, fix the matrix with diagonal entries $0$ and off-diagonal entries $1$:
		\begin{equation}
			\label{eq:J4} M \, = \, 
			\begin{bmatrix}
				\, 0 & 1 & 1 & 1 \, \\
				\, 1 & 0 & 1 & 1 \, \\
				\, 1 & 1 & 0 & 1 \,\\
				\, 1 & 1 & 1 & 0 \, \end{bmatrix}. 
		\end{equation}
		The polytrope $Q_M$ is the set of solutions to the $12$  inequalities
		$u_i - u_j \leq 1$ for $i \not= j$. It is the $3$-dimensional polytope
		shown in Figure \ref{fig:rhombic}. Namely,
		$Q_M$ is a {\em rhombic dodecahedron}, with
		$14$ vertices, $24$ edges and $12$~facets.
		The vertices  are the images in $\RR^4 / \RR {\bf 1}$ of the
		$14$ vectors~in $\{0,1\}^4 \backslash \{ {\bf 0},{\bf 1}\}$.
		Vertices  $e_i $ are blue, vertices
		$e_i + e_j$ are yellow, and vertices
		$e_i+e_j+e_k$ are~red.
		
		The order $\Lambda_M$ consists of all
		$4 \times 4$ matrices with entries in the
		valuation ring $\Ocal_K$ whose off-diagonal elements lie in the maximal ideal $\langle p \rangle$. 
		We shall see in Theorem \ref{thm:injproj}
		that the blue and red vertices encode 
		the injective modules and the projective modules of $\Lambda_M$ respectively.
		
		\begin{figure}[h]
			\vspace{-0.05cm}		
			\begin{center}
				\includegraphics[scale=0.4]{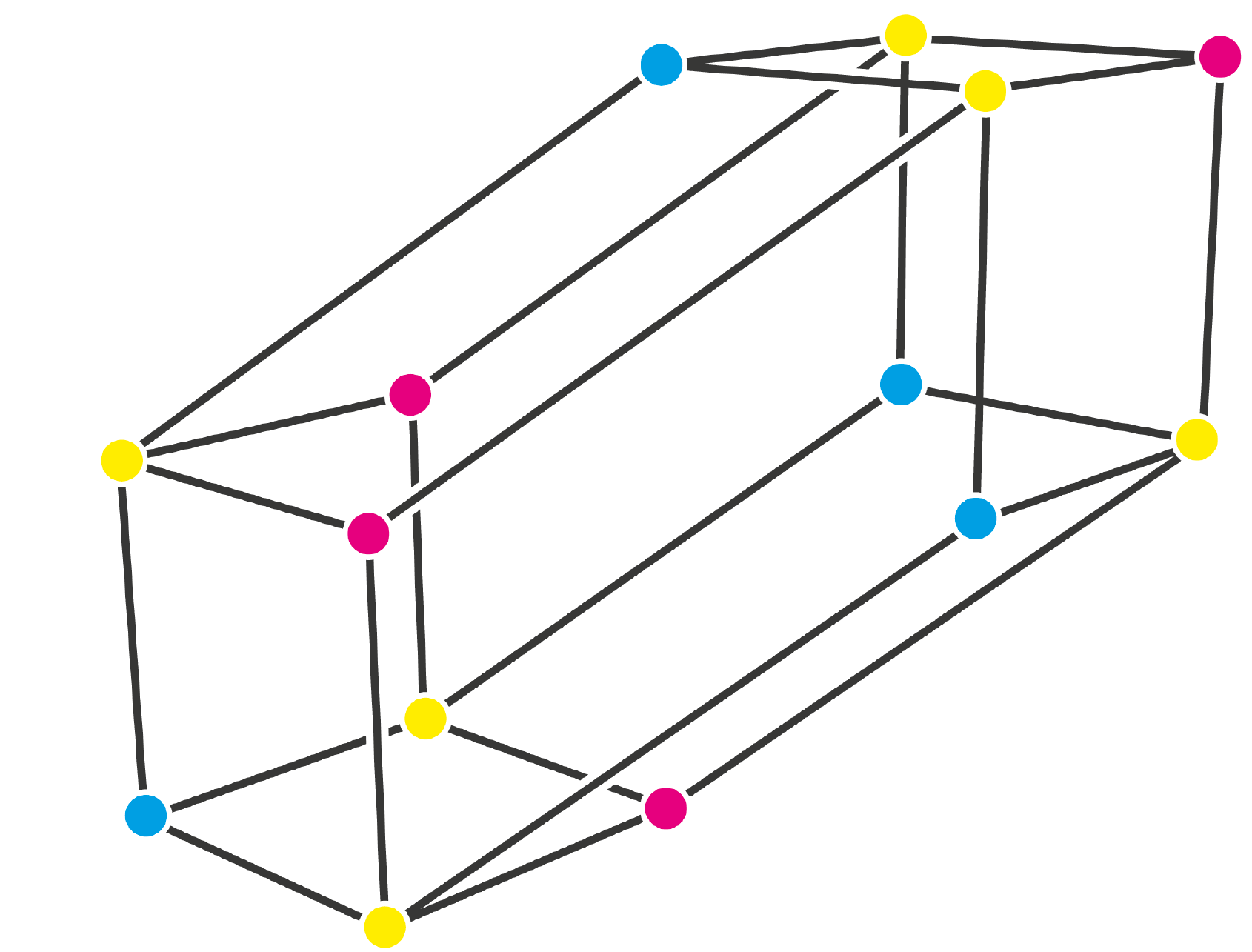} \quad
				\includegraphics[scale=0.4]{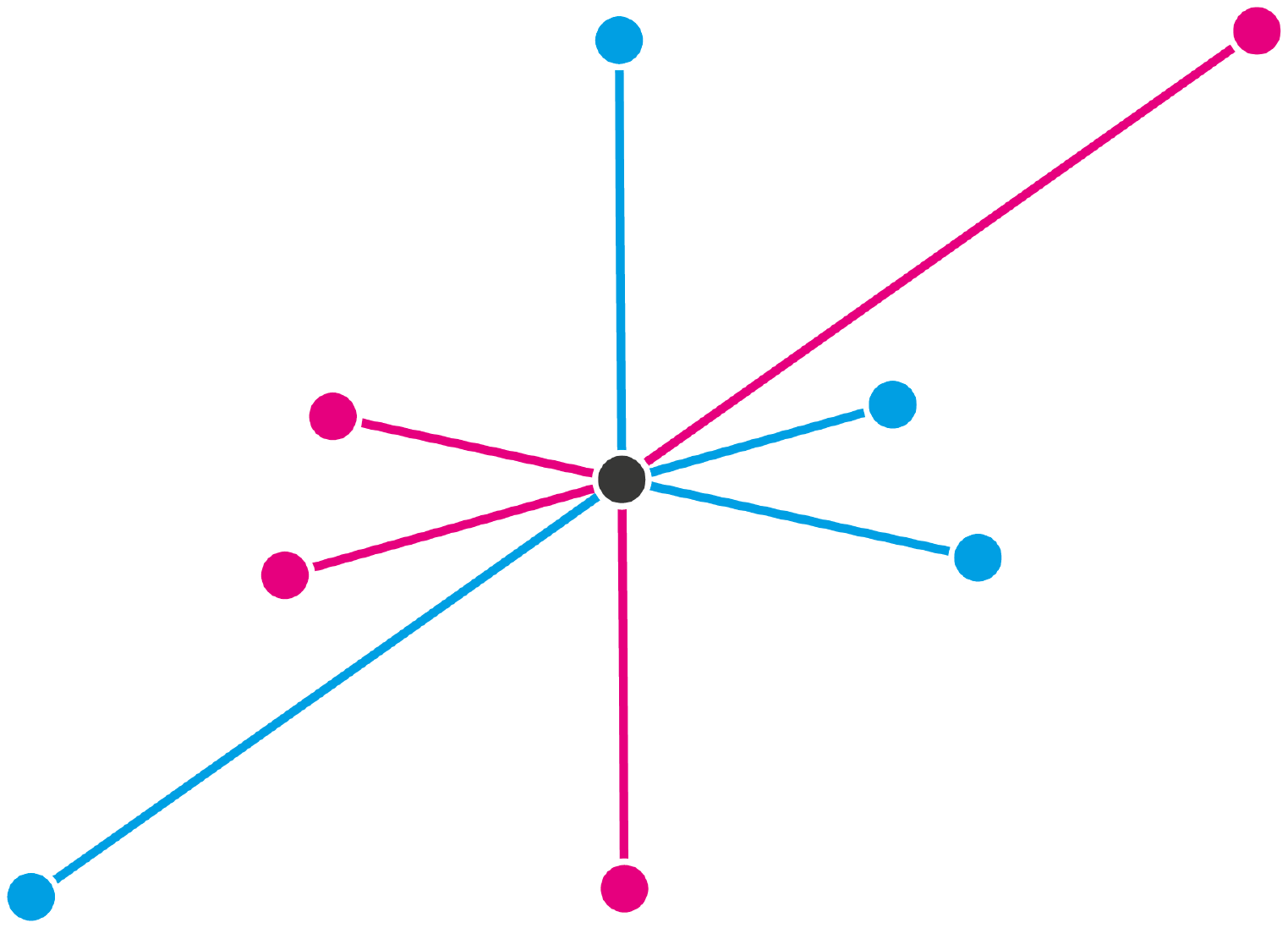}
			\end{center}
			\vspace{-0.3cm}			
			\caption{The polytrope $Q_M$ on the left is a rhombic dodecahedron.
				The four blue vertices and the four red vertices, highlighted on the right, 
				will play a special role
				for the order $\Lambda_M$.
				\label{fig:rhombic}}
		\end{figure}
	\end{example}
	
	The connection between algebra, geometry and combinatorics
	we present was pioneered by Plesken and Zassenhaus.
	Our primary source on their work is the book \cite{plesken83}.
	One objective of this article is to give an exposition of their results
	using the framework of tropical geometry \cite{joswig, maclagan}.
	But we also present a range of new results.
	Our presentation is organized as follows.
	
	Section \ref{sec2} concerns graduated orders in $K^{d \times d}$.
	In Propositions~\ref{prop:isanorder} and \ref{prop:isstable}
	we present linear inequalities that
	characterize these orders and the lattices they act on.
	These inequalities play an important role in
	tropical convexity, to be explained in Section \ref{sec3}.
	Theorem~\ref{thm:PZtropical} gives a tropical matrix formula  for  the
	Plesken-Zassenhaus order of a collection of diagonal lattices.
	
	In Section \ref{sec4} we introduce 	polytrope regions. These are convex
	cones and polyhedra whose integer points	represent graduated orders.
	Section \ref{sec5} is concerned with
	(fractional) ideals 	in an order $\Lambda_M$.
	These are parametrized by the ideal class polytrope $\mathcal{Q}_M$.
	In Section \ref{sec6} we turn to Bruhat-Tits buildings and their chambers.
	While the present study is restricted
	to Plesken-Zassenhaus orders arising from 
	one single apartment, it sets the stage for a general~theory.
    	
    Several of the results in this article were found by computations.
    The codes and all data are made available at
    	 \url{https://mathrepo.mis.mpg.de/OrdersPolytropes/index.html}.

	\section{Graduated Orders}
	\label{sec2}	
	
	By a  {\em lattice} in $K^d$ we mean a free $\Ocal_K$-submodule of rank $d$.
	Two lattices $L$ and $L'$ are equivalent if 
	$L' =  p^n L $ for some $n \in \ZZ$. We write
	$[L] = \{ p^n L : n \in \ZZ \}$ for the equivalence class of $L$.
	An \emph{order} in $K^{d\times d}$ is a lattice
	in the $d^2$-dimensional vector space $K^{d\times d}$ that is also 
	a ring. Thus, every order contains the identity matrix. 
	An order $\Lambda$  is {\em maximal} if it
	is not properly contained in any other order.
	One example of a maximal order is the matrix ring
	$$\Ocal _K^{d\times d} \,\,:=\,\, \{ X \in K^{d\times d} : 
	\val (x_{ij}) \geq 0 \mbox{ for all } 1\leq i, j \leq d \} .$$
	This is spanned as an $\Ocal _K$-lattice
	by the matrix units $E_{ij} $ where $1\leq i, j \leq d $.
	It is multiplicatively closed because $E_{ij} E_{jk} = E_{ik}$.
	We begin with some standard facts found in~\cite{plesken83}.
	The first is a natural bijection between 
	lattice classes $[L]$ in $K^d$ and maximal orders in $K^{d\times d}$.

	\begin{proposition}\label{prp:max order}
		Any order $\Lambda $ in $K^{d\times d}$ is contained
		in the endomorphism ring of a lattice $L\subset K^d$.
		The maximal orders in $K^{d\times d}$ 
		are exactly the endomorphism rings  of lattices  $L$:
		$$ \End _{\Ocal _K}(L) \,:=\, \{ X \in K^{d\times d} :  X L \subseteq L \} .$$ 			
		Two lattices $L$ and $L'$ in $K^d$ are equivalent if and only if
		$\,\End _{\Ocal _K}(L) = \End _{\Ocal _K}(L')   $.
	\end{proposition}
	
	\begin{proof}
		Let $\Lambda  = \bigoplus _{j=1}^{d^2} \Ocal _K X_j $ be
		an order in $K^{d \times d}$. 
		If we apply all the matrices $X_j$ to the {\em standard lattice}
		$L_0 =  \Ocal_K^d = \bigoplus _{i=1}^d \Ocal _K e_i $, then
		we obtain the following lattice in $K^d$:
		$$L\,\,:=\,\,\sum _{j=1}^{d^2} X_j L_0 \, \,= \,\,
		\sum_{i=1}^d \sum_{j=1}^{d^2}{\Ocal _K} X_j \,e_i .  $$
		Since $\Lambda$ is multiplicatively closed, we have
		$X_j L \subseteq L$ for all $j$.
		Therefore $\Lambda \subseteq \End _{\Ocal _K}(L)$.
		
		Endomorphism rings of lattices are orders. Indeed,  
		if $L=g L_0$ for $g \in \GL_d(K)$, then
		\begin{equation}\label{eq:EndOKL}
			\End_{\Ocal _K}(L) \,=\, g \End _{\Ocal _K}(L_0) g^{-1} 
			\,		= \,\,g \,\Ocal _K^{d\times d} g^{-1}.
		\end{equation}
		This is	 a ring, and it is spanned as an $\Ocal _K$-lattice by 
		$\{ gE_{ij}g^{-1} : 1\leq i,j\leq d \}$.
		This allows to conclude that the maximal orders 
		are exactly the endomorphism rings of lattices. 
	\end{proof}
	
	\smallskip
	
	For $u \in \ZZ^d$ we abbreviate $\, g_u = {\rm diag}( p^{u_1},p^{u_2}, \ldots,p^{u_d})$.
	This diagonal matrix transforms the standard lattice $\Ocal_K^d$ to $L_u = g_u \Ocal_K^d$. The endomorphism ring
	${\rm End}_{\Ocal_K}(L_u)$ is the maximal order in (\ref{eq:EndOKL}).
	Let $M(u) $ denote the $d \times d$ matrix whose
	entry in position $(i,j)$ equals $u_i - u_j$.
	
	\begin{lemma} The endomorphism ring of the  lattice
		$L_u$ is given by valuation inequalities:
		\begin{equation}
			\label{eq:EndOval}
			{\rm End}_{\Ocal_K}(L_u) \,\, = \,\, \Lambda_{M(u)} \,\, = \,\,
			\{\, X \in K^{d \times d} \,: \, {\rm val}(X) \geq M(u) \,\} . 
		\end{equation}
	\end{lemma}
	
	\begin{proof}
		The elements of  ${\rm End}_{\Ocal_K}(L_u)$
		are the matrices $X = g_u Y g_u^{-1}$ where $Y \in \Ocal_K^{d \times d}$.
		Writing $X = (x_{ij})$ and $Y = (y_{ij})$, the equation
		$X = g_u Y g_u^{-1}$  means that $x_{ij} = p^{u_i-u_j} y_{ij} $ for all $i,j$.
		The condition ${\rm val}(y_{ij}) \geq 0$ is equivalent to
		${\rm val}(x_{ij}) \geq u_i - u_j$. Taking the conjunction over all $(i,j)$, we conclude 
		that ${\rm val}(Y) \geq 0$ is equivalent to the desired inequality ${\rm val}(X) \geq M(u)$.
	\end{proof}
	
	The matrices $M(u)$ are characterized by the following two properties.
	All diagonal entries are zero and the tropical rank
	is one, cf.~\cite[Section 5.3]{maclagan}. 
	What happens if we replace $M(u)$ 
	in (\ref{eq:EndOval})
	by an arbitrary matrix
	$M\in \ZZ^{d \times d}$?
	Then we get the set
	$\Lambda_M$ from the Introduction.
	
	\begin{remark} For any matrix $M \in \ZZ^{d \times d}$,
		the set $\Lambda_M$ is a lattice in $K^{d \times d}$. It is generated
		as an $\Ocal_K$-module by the matrices
		$ p^{m_{ij}} E_{ij} $ for $1 \leq i,j \leq d$.
		The lattice $\Lambda_M$ may not be an~order.
	\end{remark}
	
	Write  $\ZZ^{d \times d}_0$ for the set of integer matrices $M$
	with zeros on the diagonal, i.e.~$m_{ii} = 0$ for all~$i$. If $M$ 
	lies in $\ZZ^{d \times d}_0 $
	then $\Lambda_M$ contains the identity matrix, but may still not be an~order.

	\begin{example}
		Let $K = \QQ$ with the $p$-adic valuation, for some prime $p \geq 5$.
		For $d=3$, set 
		$$ \begin{small} M = \begin{bmatrix} 0 & 0  & 1  \\
				0 & 0 & 0\\
				0 & 0 & 0  \end{bmatrix} \quad {\rm and} \quad
			X = \begin{bmatrix} 1 & 1 & p \\ 1 & 1 & 1 \\ 1 & 1 & 1
			\end{bmatrix} , \quad {\rm so} \quad
			X^2 = \begin{bmatrix} 2{+}p & 2{+}p & 1 {+} 2p \\ 3 & 3 & 2{+}p \\ 3 & 3 & 2{+}p \end{bmatrix}.  
		\end{small} $$
		Since ${\rm val}(X) = M$ and ${\rm val}(X^2) = 0$,  we have $X \in \Lambda_M$ but
		$X^2 \not\in \Lambda_M$. So
		$\Lambda_M$ is not an~order.
	\end{example}
	
	The inequalities derived in the next two propositions are the main points of this section.
	These results are due to Plesken \cite{plesken83}. He states them in
	\cite[Definition~II.2]{plesken83} and  \cite[Remark~II.4]{plesken83}.
	The orders $\Lambda_M$ in Proposition \ref{prop:isanorder}
	are called {\em graduated orders} in \cite{plesken83}. 
	They are also known as {\em tiled orders} \cite{DKKP, tiled}, {\em split orders} \cite{shemanske} or {\em monomial orders} \cite{monomial}.
	A graduated order $\Lambda_M$ is
	in {\em standard form} if $M \geq 0$ and $m_{ij} + m_{ji} > 0$
	for $i \not= j$.  		
	
	\begin{proposition} 
		\label{prop:isanorder}
		Given $M = (m_{ij})$  in $\ZZ^{d \times d}_0$, the
		lattice $\Lambda_M$ is an order in $K^{d \times d}$ if and only if
		\begin{equation}
			\label{eq:polytrope1}
			m_{ij} + m_{jk} \geq m_{ik} \quad
			\hbox{for all} \,\, \,1 \leq i,j,k \leq d. 
		\end{equation}
	\end{proposition}
	
	\begin{proof}
		To prove the if direction, we assume (\ref{eq:polytrope1}).
		Our hypothesis $m_{ii} = 0$ ensures that $\Lambda_M$ contains the identity matrix, so
		$\Lambda_M$ has a multiplicative unit.
		Suppose $X,Y \in \Lambda_M$. Then the $(i,k)$ entry of $XY$ equals
		$\sum_{j=1}^d x_{ij} y_{jk}$. This is a scalar in $K$ whose valuation is
		at least $m_{ij} + m_{jk}$ for some index $j$.
		Hence it is greater than or equal to $m_{ik}$ since (\ref{eq:polytrope1}) holds.
		
		For the only-if direction, suppose $m_{ij} + m_{jk} < m_{ik}$.
		Then $X = p^{m_{ij}} E_{ij}$  and 
		$Y = p^{m_{jk}} E_{jk}$  are in $\Lambda_M$.
		However, $XY = p^{m_{ij}+m_{jk}} E_{ik}$  is not in $\Lambda_M$
		because its entry in position $(i,k)$ has valuation less than $m_{ik}$.
		Hence $\Lambda_M$ is not multiplicatively closed, so it is not an order.
	\end{proof}
	
	Fix  $M$ that satisfies (\ref{eq:polytrope1}). The graduated order
	$\Lambda_M$ is an $\Ocal_K$-subalgebra 
	of $K^{d\times d}$. 
	It is therefore natural to ask which lattices in $K^d$ are $\Lambda_M$-stable.
	
	\begin{proposition}\label{prop:isstable}
		A lattice $L$ is stable under $\Lambda_M$ if and only if
		$L = L_u$ with $u \in \ZZ^d$ satisfying
		\begin{equation}
			\label{eq:polytrope2} u_i - u_j \,\leq\, m_{ij}  \quad {\rm for} \quad 1 \leq i,j \leq d . 
		\end{equation}
		Moreover, if $u,u' \in \ZZ^d$ satisfy (\ref{eq:polytrope2}), then 
		the diagonal lattices $L_u$ and $L_{u'}$  are isomorphic as $\Lambda_M$-modules if and only if
		they are equivalent, i.e.~$u = u'$ in the 
		quotient space $\RR^d/\RR {\bf 1}$.
	\end{proposition}
	
	\begin{proof} 
		Fix a lattice $L$ and let $u = (u_1,\ldots,u_d)$ be defined by 
		$u_i = {\rm min}\{ {\rm val}(b_i) : b \in L \}$.
		Then $L \subseteq L_u$ because every $b \in L$ is an 
		$\Ocal_K$-linear combination of the standard basis of $L_u$,
		namely $\,b = \sum_{i=1}^d b_i e_i = \sum_{i=1}^d (b_i \,p^{-u_i}) \,p^{u_i} e_i $.
		Suppose that $L$ is $\Lambda_M$-stable. Since $m_{ii} = 0$,
		we have $E_{ii} \in \Lambda_M$. Hence $E_{ii}\, b = b_i e_i \in
		L$ for every $ b \in L$. This implies $L_u \subseteq L$ and hence~$L = L_u$.				
		Applying $p^{m_{ij} } E_{ij} \in \Lambda_M$ to $p^{u_j} e_j \in L_u$, we see that
		$p^{m_{ij} + u_j} e_i$ lies in $L_u$, and this implies $m_{ij} + u_j \geq u_i$.
		Hence (\ref{eq:polytrope2}) holds. Conversely, suppose that (\ref{eq:polytrope2}) holds.
		Then the generator $p^{m_{ij}} E_{ij}$
		of $\Lambda_M$ maps each basis vector  $p^{u_k} e_k$ of $L_u$ either to zero
		(if $j \not= k)$, or to 
		$p^{m_{ik}+u_k}e_i\in L_u$.
		This proves the first assertion.
		
		For the second assertion, let $u,u' \in \ZZ^d$ satisfy (\ref{eq:polytrope2}).
		Since multiplication by  $\alpha \in K^*$  is an isomorphism of $\Ocal_K$-modules,
		the if-direction is clear. Conversely, if
		$L_u$ and $L_{u'}$ are isomorphic, then there exists $g \in {\rm GL}_d(K)$ such that
		$L_{u'} = g L_u$ and $g X = X g$ for all $X \in \Lambda_M$.
		Pick $s \in \ZZ_{>0}$ such that $p^s \Ocal_K^{d \times d} \subset \Lambda_M$.
		Then $g$ commutes with every matrix in $p^s \Ocal_K^{d \times d}$.
		This implies that $g$ is central in $\Ocal_K^{d \times d}$, and therefore
		$g$ is a multiple of the identity matrix.
	\end{proof}

	\section{Bi-tropical Convexity}
	\label{sec3}	
	
	We now develop the relationship between graduated orders and		
	tropical mathematics \cite{joswig,
		maclagan}. Both
	the {\em min-plus algebra}
	$ (\,\RR, \minplus  ,\odot)$
	and the {\em max-plus algebra}
	$ (\,\RR, \maxplus ,\odot)$ will be used. Its
	arithmetic operations are the minimum, maximum, and classical addition 
	of real numbers:
	$$  a \minplus b \, = \, {\rm min}(a,b) \, , \,\,\,
	a \maxplus b \, = \, {\rm max}(a,b) \, , \,\,\,
	a \odot b\, = \, a + b \quad\,\, {\rm for} \quad a,b \in \RR. $$
	If $M$ and $N$ are real matrices, and the number of columns of $M$ equals
	the number of rows of $N$, then
	we write $M \mintimes N$ and $M \maxtimes N$ for their
	respective matrix products in these algebras.
	
	\begin{example} \label{ex:22MN} Consider the $2 \times 2$ matrices
		$\begin{small}  M = \begin{bmatrix} 0 & 1 \\ 2 & 0  \end{bmatrix}   \end{small}$
		and  $\begin{small}  N  = \begin{bmatrix} 1 & 0 \\ 0 & 0 \end{bmatrix} \end{small} $.
		We find that
		\vspace{-0.2cm}
		$$ \begin{matrix} 
			M \mintimes M =   \begin{bmatrix}  
				0 & 1 \\
				2 & 0  \end{bmatrix} \, , & 
			M \mintimes N =   \begin{bmatrix}  
				1 & 0 \\
				0 & 0  \end{bmatrix} \, , & 
			N \mintimes M =   \begin{bmatrix}  
				1 & 0 \\
				0 & 0  \end{bmatrix} \, , & 
			N \mintimes N =   \begin{bmatrix}  
				0 & 0 \\
				0 & 0  \end{bmatrix} \, , \medskip \\
			M \maxtimes M =   \begin{bmatrix}  
				3 & 1 \\
				2 & 3  \end{bmatrix} \, , & 
			M \maxtimes N =   \begin{bmatrix}  
				1 & 1 \\
				3 & 2  \end{bmatrix}  \, ,& 
			N \maxtimes M =   \begin{bmatrix}  
				2 & 2 \\
				2 & 1  \end{bmatrix}  \, ,& 
			N \maxtimes N =   \begin{bmatrix}  
				2 & 1 \\
				1 & 0  \end{bmatrix}  .
		\end{matrix}
		$$
	\end{example}
	There are two flavors of tropical convexity \cite[Section 5.2]{maclagan}.
	A subset of $\RR^d$ is {\em min-convex}  if it is closed under linear combinations in the min-plus algebra,
	and {\em max-convex}  if the same holds for the max-plus algebra.
	Thus convex sets are images of matrices under linear maps.
	
	We are especially interested in bi-tropical convexity 
	in the
	ambient space $\RR^d/\RR {\bf 1}$. This 
	is ubiquitous in \cite[Section 5.4]{joswig} and \cite{maclagan}. Joswig  \cite[Section~1.4]{joswig} calls it
	the {\em  tropical projective torus}.
	At a later stage, we also work in the corresponding matrix space
	$\RR^{d \times d}/ \RR {\bf 1}$.
	
	Let $\RR^{d \times d}_0$ denote the space of real 
	$d \times d$ matrices with zeros on the diagonal, which is a real $(d^2-d)$-dimensional vector space with lattice $\ZZ^{d \times d}_0$.
	For $M= (m_{ij}) $ in $\RR^{d \times d}_0$, we define 
	\begin{equation}
		\label{eq:polytrope}
		Q_M  \,\, = \,\, \bigl\{ u \,\in \RR^d/\RR{\bf 1} \,:\, u_i - u_j \,\, \leq \, m_{ij} \,\,\,
		\hbox{for} \,\, 1 \leq i,j \leq d \,\bigr\}. 
	\end{equation}
	Such a set is known as a {\em polytrope} in tropical geometry \cite{joskul, maclagan}.
	Other communities use the terms
	{\em alcoved polytope} and {\em weighted digraph polytope}.
	We note that  $Q_M$ is both
	min-convex and max-convex \cite[Proposition 5.30]{joswig}
	and, being a polytope, it is also classically convex.
	
	Using tropical arithmetic, the linear inequalities in
	(\ref{eq:polytrope1}) can be written concisely as follows:
	\begin{equation}
		\label{eq:polytrope3} M \mintimes M \,= \, M. 
	\end{equation}
	Thus, $M$ is {\em min-plus idempotent}. This holds for $M$ in Example \ref{ex:22MN}.
	Joswig's book \cite[Section 3.3]{joswig} uses the term {\em Kleene star} for matrices
	$M \in \RR^{d \times d}_0$ with (\ref{eq:polytrope3}).
	Propositions \ref{prop:isanorder} and \ref{prop:isstable} imply:
	
	\begin{corollary} The lattice $\Lambda_M$ is an order
		in $K^{d \times d}$ if and only if (\ref{eq:polytrope3}) holds.
		In this case, the integer points $u$ in the polytrope $Q_M$
		are in bijection with  the isomorphism classes of $\Lambda_M$-lattices $L_u$.
		Here, by a {\em $\Lambda_M$-lattice} we mean a 
		$\Lambda_M$-module that is also a lattice in $K^d$.
	\end{corollary}
	
	Let $\Gamma = \{ L_1, \ldots, L_n \}$ be a finite set of lattices in $K^d$, which might be taken up to
	equivalence.  The intersection of two orders in $K^{d \times d}$ is again an order.
	Hence the  intersection
	\begin{equation}
		\label{eq:PZ1} {\rm PZ}(\Gamma) \,\, = \,\,
		{\rm End}_{\Ocal_K}(L_1) \,\cap \,\cdots \, \cap \,{\rm End}_{\Ocal_K}(L_n)
	\end{equation}
	is an order in $K^{d \times d}$. We call $\PZ(\Gamma)$ the
	{\em Plesken-Zassenhaus order} 
	of the configuration $\Gamma$.
	
	In the following we assume that each $L_i$ is
	a {\em diagonal lattice}, i.e.~$L_i = L_{u^{(i)}}$ for $u^{(i)} \in \ZZ^d$.
	Our next result involves
	a curious mix of max-plus
	algebra and min-plus algebra.
	
	\begin{theorem} \label{thm:PZtropical}
		Let  $\,\Gamma = \{L_{u^{(1)}}, \ldots, L_{u^{(n)}}\}\,$ be any configuration of
		diagonal lattices in $K^d$. Then its Plesken-Zassenhaus order 
		$PZ(\Gamma)$ coincides with the
		graduated order $\Lambda_M$ where
		\begin{equation}
			\label{eq:maxplusM}
			M \,\,\, = \,\,\, M(u^{(1)}) \,\maxplus \,
			M(u^{(2)})  \,\maxplus \,\cdots\, \maxplus \, M(u^{(n)}) .
		\end{equation}
		This max-plus sum of tropical rank one
		matrices is min-plus idempotent, i.e.~(\ref{eq:polytrope1}) and (\ref{eq:polytrope3})~hold.	\end{theorem}
	
	\begin{proof}
		We regard $\Gamma$ as a configuration in $\RR^d/\RR {\bf 1}$.
		By construction, $M$ is the entrywise smallest matrix such that
		$\Gamma $ is contained in the polytrope $Q_M$.
		From \cite[Lemma~3.25]{joswig} the matrix $M$ is a Kleene star, that is (\ref{eq:polytrope1}) and (\ref{eq:polytrope3}) hold. 
		The intersection in
		(\ref{eq:PZ1}) is defined by the conjunction of the $n$
		inequalities ${\rm val}(X) \geq M(u^{(i)})$, which is equivalent to ${\rm val}(X) \geq M$.
	\end{proof}
	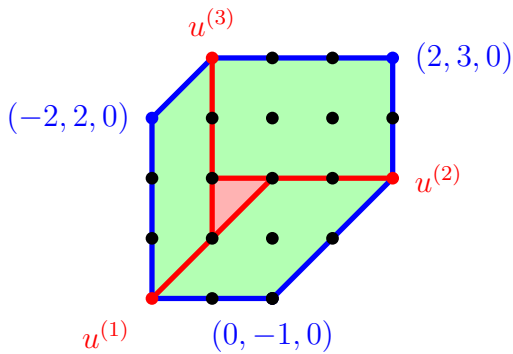
\begin{figure}[h]		
		\vspace{-0.22cm}
		\begin{center}
			\begin{tikzpicture}[scale=0.8]
				\fill[green!30]  (-2,-5) -- (0,-5) -- (2,-3) -- (2,-1) -- (-1,-1) -- (-2, -2) -- (-2 , -5);
				\fill[red!30]  (-1,-3) -- (0,-3) -- (-1,-4) -- (-1,-3);
				\draw[thick, red, line width=0.7mm] (-2, -5) -- (0, -3);
				\draw[thick, red, line width=0.7mm] (-1,-3)  -- (2,-3);
				\draw[thick, red, line width=0.7mm] (-1, -4) -- (-1, -1);
				\draw[thick, blue, line width=0.7mm] (-2,-5) -- (0,-5) -- (2,-3) -- (2,-1) -- (-1,-1) -- (-2, -2) -- (-2 , -5);						
				\fill[red] (-2,-5) circle (3pt) node[below left= 1.5mm]{$u^{(1)}$};
				\fill[red] (2,-3) circle (3pt)  node[right = 1.5mm]{$u^{(2)}$};
				\fill[red] (-1,-1) circle (3pt) node[above = 1.5mm]{$u^{(3)}$}; 
				\fill[blue] (2,-1) circle (3pt) node[right = 1.5mm]{$(2,3,0)$};
				\fill[blue] (0,-5) circle (3pt) node[below = 1.5mm]{$(0,-1,0)$}; 		
				\fill[blue] (-2,-2) circle (3pt) node[left = 1.5mm]{$(-2,2,0)$}; 
				\fill[black] (-1,-2) circle (3pt); 
				\fill[black] (-1,-3) circle (3pt); 
				\fill[black] (-1,-4) circle (3pt); 
				\fill[black] (-1,-5) circle (3pt); 
				\fill[black] (-2,-4) circle (3pt); 
				\fill[black] (-2,-3) circle (3pt); 
				\fill[black] (0,-1) circle (3pt); 
				\fill[black] (0,-2) circle (3pt); 
				\fill[black] (0,-3) circle (3pt); 
				\fill[black] (0,-4) circle (3pt); 
				\fill[black] (0,-5) circle (3pt); 
				\fill[black] (1,-1) circle (3pt); 
				\fill[black] (1,-2) circle (3pt); 
				\fill[black] (1,-3) circle (3pt); 
				\fill[black] (1,-4) circle (3pt); 
				\fill[black] (0,-5) circle (3pt); 
				\fill[black] (2,-2) circle (3pt); 
			\end{tikzpicture}
		\end{center}
		\vspace{-0.4cm}		
		\caption{\label{fig:greenhexagon} A polytrope with
			three min-plus vertices (blue) and three
			max-plus vertices (red).}
	\end{figure}

	\begin{example}\label{ex:hexagon}
		For $d=n=3$, fix $u^{(1)} = (-2 , -1,0) $, $
		u^{(2)} = (2 , 1,0) $, $u^{(3)} = (-1 , 3,0)\,$
		in $\RR^3/\RR {\bf 1}$. The configuration $\Gamma = \{ u^{(1)},
		u^{(2)},u^{(3)}\}$ consists of the three red points in  Figure \ref{fig:greenhexagon}.			
		The red diagram is their min-plus convex hull. This tropical triangle
		consists of a classical triangle together with three red line segments connected to $\Gamma$. 
		This red min-plus triangle is not convex. The green shaded hexagon is the
		polytrope spanned by $\Gamma$. 
		By \cite[Remark 5.33]{joswig}, this is the geodesic convex hull of~$\Gamma$.
		It equals $Q_M$ where $M$ is computed by 
		(\ref{eq:maxplusM}):
		$$ M \,\,\,=\,\,\, (u^{(1)})^t \odot (-u^{(1)}) \,\,\,\maxplus \,\,\,
		(u^{(2)})^t \odot (-u^{(2)}) \,\,\,\maxplus \,\,\,
		(u^{(3)})^t \odot (-u^{(3)}) \,\,\,=\,\,\, \begin{small}
			\begin{bmatrix}
				0 & 1 & 2 \\
				4 & 0 & 3 \\
				2 & 1 & 0 
			\end{bmatrix}. \end{small}
		$$
		The polytrope $Q_M$ is both a min-plus triangle and a max-plus triangle.
		Its min-plus vertices, shown in blue, are
		equal in $\RR^3/\RR {\bf 1}$ to the columns of $M$.
		Its max-plus vertices, shown in red, are the points
		$u^{(i)}$. These are equal in 
		$\RR^3/\RR {\bf 1}$ to the columns of $-M^t$; cf.~Theorem \ref{thm:injproj}. Moreover, the three green cells correspond to the collection of homothety classes of  lattices contained in $u^{(i)}\minplus u^{(j)}$ and containing $u^{(i)}\maxplus u^{(j)}$, for each choice of $i\neq j$.
	\end{example}

	\begin{remark} \label{rem:indecomposable}
		All lattices $L_u$ for $u \in Q_M$ are indecomposable as $\Lambda_M$-modules, cf.~\cite{plesken83}.
		This is no longer true if $\RR$ is enlarged to
		the tropical numbers $\RR \cup \{ \infty \}$.
		The combinatorial theory of polytropes in \cite{joswig}
		is set up for this extension, and it indeed makes sense
		to study orders~$\Lambda_M$ with $m_{ij} = \infty$.
		While we do not pursue this here, our approach
		would extend to that setting.
	\end{remark}

	\begin{example} Set $d=4$. 
		The rhombic dodecahedron in Example \ref{ex:rhombic}
		was called the {\em pyrope} in
		\cite[Figure 4]{joskul} and can be seen as a unit ball with respect to the tropical metric, cf.~\cite[\S 3.3]{CGQ04}. This $Q_M$
		is a tropical tetrahedron for both min-convexity and max-convexity.
		The respective vertices are shown in red and blue in Figure~\ref{fig:rhombic}.
		We have $\Lambda_M = {\rm PZ}(\Gamma)$
		where $\Gamma$ is either set of four vertices.
		The $\Lambda_M$-lattices $L_u$ correspond to the
		$15$ integer points in~$Q_M$.
	\end{example}

	\section{Polytrope Regions}
	\label{sec4}
	
	We next introduce a cone that parametrizes all graduated orders $\Lambda_M$.
	Following Tran~\cite{tran}, the {\em polytrope region} $\mathcal{P}_d  $ is  the
	set of all min-plus idempotent matrices $M \in \RR^{d \times d}_0$. Thus, $\mathcal{P}_d$ is the
	$(d^2-d)$-dimensional convex polyhedral cone defined by the linear inequalities
	in (\ref{eq:polytrope1}). The equations $m_{ik} = m_{ij} + m_{jk}$ define the
	cycle space of the complete bidirected graph~$\mathcal{K}_d$. This is   the lineality space of $\mathcal{P}_d$.
	Modulo this $(d-1)$-dimensional space, the polytrope region~$\mathcal{P}_d$ is a pointed cone of 
	dimension $(d{-}1)^2$. We view it as a polytope of dimension $d^2-2d$. 
	Each inequality $m_{ik} \leq m_{ij} + m_{jk}$ is facet-defining, so the
	number of facets of $\mathcal{P}_d$ is $d(d-1)(d-2)$.
	
	It is interesting but difficult to list the vertices of $\mathcal{P}_d$ and to 
	explore the face lattice. The same problem was studied in \cite{avis}
	for the {\em metric cone}, which is the restriction of $\mathcal{P}_d$ 
	to the subspace of symmetric matrices in $\RR^{d \times d}_0$.
	A website maintained
	by Antoine Deza \cite{deza} reports that the number of rays of the metric cone
	equals  $3,7, 25, 296, 55226, 119269588$ for $d=3,4,5,6,7,8$.
	We here initiate the census for the polytrope region.
	The following tables report the
	size of the orbit, the number of incident facets, and a representative  matrix $[m_{ij}]$.
	Here orbit and representatives refer to the natural action of the symmetric group $S_d$
	on $\mathcal{P}_d$.
	The matrices $[m_{ij}]$ in $\ZZ^{3 \times 3}_0$ are written in the vectorized format
	$[m_{12}m_{13}m_{21}m_{23}m_{31}m_{32}]$.
	
	\begin{proposition}
		The  polytope $\mathcal{P}_3$ is a bypramid,
		with f-vector $(5,9,6)$. Its five vertices are
		\[
		3 \,,\, 4\  [001100]\ \    \textup{ and }\ \  2\, ,\, 3 \  [001110].
		\]
		The polytope $\mathcal{P}_4$ has the f-vector
		$(37,327,1140,1902,1680,808,204,24)$. Its $37$ vertices are
		\[ \!\!\! \begin{small} \begin{matrix} 
				12,  10 \!\!\! &   [1 1 1 0 1 1 0 0 1 0 0 1] & \quad
				6 , 12 \!\!\! &   [1 1 1 0 1 1 0 0 1 0 0 0]  & \quad
				12 , 14 \!\!\! &   [0 1 1 0 1 1 0 0 1 0 0 0] \\
				3 , 16 \!\!\! &    [0 1 1 0 1 1 0 0 0 0 0 0]   & \quad
				4  ,18 \!\!\! &    [1 1 1 0 0 0 0 0 0 0 0 0]   .
		\end{matrix} \end{small}
		\]
		The corresponding polytropes $Q_M$ are pyramid, tetrahedron, triangle, segment, and segment.		
		The $15$-dimensional polytope $\mathcal{P}_5$ has
		$2333$ vertices in $33$ symmetry classes. These classes are
		\[ \!\!\! \begin{tiny} \begin{matrix}
				5 , 48 \!\!\! & [00000000000000001111] & 
				10 , 18 \!\!\! & [00001001211121111100] & 
				10 , 42 \!\!\! & [00000000000011101110] \\
				20 , 15 \!\!\! & [00002012323231012201] &
				20 , 21 \!\!\! & [00001000110021112111] &
				20 , 39 \!\!\! & [00000000000011101111] \\
				24 , 20 \!\!\! & [00001001210122111110] & 
				24 , 30 \!\!\! & [00001000110011101111] &
				30 , 24 \!\!\! & [00001000110121111110] \\
				30 , 30 \!\!\! & [00000000110011111111] &
				30 , 30 \!\!\! & [00000000110111111110] &
				30 , 36 \!\!\! & [00000000110011001111] \\
				40 , 18 \!\!\! & [00002000221222212212] &
				60 , 18 \!\!\! & [00001000210122112110] &
				60 , 18 \!\!\! & [00001001210122121100] \\
				60 , 22 \!\!\! & [00001000110122111110] &
				60 , 27 \!\!\! & [00001000110011102111] &
				60 , 29 \!\!\! & [00000000110011102211] \\
				60 , 33 \!\!\! & [00000000110011101111] &
				120 , 16 \!\!\! & [00001001220132122110] &
				120 , 17 \!\!\! & [00001001210122122110] \\
				120 , 18 \!\!\! & [00001001210122112110] &
				120 , 18 \!\!\! & [00001001210122122210] &
				120 , 18 \!\!\! & [00001001210222122110] \\
				120 , 18 \!\!\! & [00001001220132213210] &
				120 , 19 \!\!\! & [00001000210022103221] &
				120 , 19 \!\!\! & [00001000210122122110] \\
				120 , 19 \!\!\! & [00001001210122212210] &
				120 , 22 \!\!\! & [00001000110021102221] &
				120 , 22 \!\!\! & [00001000110122121110] \\
				120 , 23 \!\!\! & [00001000110021102211] &
				120 , 23 \!\!\! & [00001000110021102222] &
				120 , 25 \!\!\! & [00001000110011102211]
			\end{matrix}
		\end{tiny}
		\]
	\end{proposition}
	
	\begin{proof}
		This was found by computations with {\tt Polymake} \cite{polymake}; see our {\tt mathrepo} site.
	\end{proof}
	
	\begin{remark}
		The integer matrices $M$ in the 	polytrope region $\mathcal{P}_d$ represent the graduated 
		orders $\Lambda_M \subset K^{d \times d}$. The data above enables us to sample from these orders. A variant of $\mathcal{P}_d$ that assumes nonnegativity constraints was studied in \cite{DDP}, which offers additional data. We also refer to \cite{DKKP} for a study of the cone of polytropes from the perspective of semiring theory. 	

\end{remark}	
	
	Our next result relates the structure of a polytrope $Q_M$ to that of its graduated order~$\Lambda_M$.
	
	\begin{theorem} \label{thm:injproj}
		Let $M \in \mathcal{P}_d$ be in standard form.
		The		$(d-1)$-dimensional
		polytrope $Q_M$ is both a min-plus simplex
		and a max-plus simplex. 
		The min-plus vertices $u$ are the columns of $M$.
		They represent precisely those modules $L_u$ over the order $\Lambda_M$ that are projective.
		The max-plus vertices $v$ are the columns of $-M^t$,
		and they represent the injective $\Lambda_M$-modules~$L_v$.
	\end{theorem}
	
	\begin{proof}
		Thanks to \cite[Theorem~7]{joskul}, full-dimensional polytropes are tropical simplices,
		with vertices given by the columns of the defining matrix $M$.
		We know from  bi-tropical convexity \cite[Proposition 5.30]{joswig} that
		$Q_M$ is both min-convex and max-convex, so it is a simplex in both ways.
		This duality corresponds to swapping $M$ with its negative transpose $-M^t$.
		Note its appearence in \cite[Theorem 5.2.21]{maclagan}.
		The connection to projective and injective modules appears in
		parts (v) and (vii)  of \cite[Remark~II.4]{plesken83}. For completeness, we sketch a proof.
		
		Recall that a module is projective if and only if it is a direct summand of a free module.
		Let $m^{(1)},\ldots,m^{(d)}$ denote the columns of $M$.  The
		lattice associated to the $j$-th column equals
		$$ L_{m^{(j)}} \,\,\,= \,\,\, \bigl\{\,x \in K^d \,:\, {\rm val}(x_i) \geq  m_{ij} \,\,\,
		\hbox{for}\,\, i = 1,\ldots,d \,\bigr\}. $$
		Taking the direct sum of these $d$ lattices gives the following identification of $\Ocal_K$-modules:
		\begin{equation}
			\label{eq:decomposed}
			\Lambda_M \,\, = \,\,  L_{m^{(1)}} \,\oplus  \,L_{m^{(2)}} \,\oplus \,\cdots \,\oplus \,L_{m^{(d)}} .
		\end{equation}
		We see that $L_{m^{(j)}}$ is a direct summand of the free rank one module $\Lambda_M$, so it
		is projective.
		
		Conversely, let $P$ be any indecomposable projective $\Lambda_M$-module.
		Then $P \oplus Q \cong \Lambda_M^r$ for some module $Q$
		and some $r \in \ZZ_{>0}$. The module $\Lambda_M^r$ decomposes
		into $r\cdot d$ indecomposables, found by aggregating $r$ copies of
		(\ref{eq:decomposed}). By the Krull-Schmidt Theorem,
		such decompositions are unique up to isomorphism, and hence $P$
		is isomorphic to $L_{m^{(j)}}$ for some~$j$.
		
		A $\Lambda_M$-module $P$ is projective if and only if ${\rm Hom}_{\Ocal_K}(P,\Ocal_K)$
		is an injective $\Lambda_M$-module, but now with the action on the right.
		The decomposition (\ref{eq:decomposed}) dualizes gracefully. 
		We derive the assertion  for injective modules
		by similarly dualizing all steps in the argument above.
	\end{proof}
	
In relation to Theorem~\ref{thm:injproj} we remark that the columns and negative rows of $M$ also have a natural interpretation as potentials in combinatorial optimization; cf.\ \cite[Theorem~3.26]{joswig}.

	\begin{example} \label{ex:rhombic2}
		The columns of the matrix $M$  in  Example \ref{ex:rhombic} are
		the negated unit vectors $-e_i$. The columns of $-M^t$ are the
		unit vectors $e_i$.
		Our color coding in Figure \ref{fig:rhombic}
		exhibits the two structures of  $Q_M$ as a tropical tetrahedron in $\RR^4/\RR {\bf 1}$.
		The four red points are the min-plus vertices, giving the
		projective $\Lambda_M$-modules.
		The four blue points are the max-plus vertices.
	\end{example}
	
	Given any min-plus idempotent matrix $M  \in \mathcal{P}_d$, 
	we define its {\em truncated polytrope region}
	\begin{equation}
		\label{eq:TPR}
		\mathcal{P}_d(M) \, = \, \{ N \in \mathcal{P}_d \,: \, N \leq M \}. 
	\end{equation}
	This polytope has dimension $d^2-d$ if
	$M$ is in the interior of $\mathcal{P}_d$.
	It parametrizes all
	subpolytropes of $Q_M$, i.e.\ all the
	polytropes $Q_N$  contained in
	$Q_M$, as the following lemma shows.
	
	\begin{lemma}
		Given matrices $M$ in $\mathcal{P}_d$ and 
		$N$ in $\RR_0^{d\times d}$ such that $Q_N\subseteq Q_M$, there exists 
		a matrix $C$ in the truncated polytrope region $\mathcal{P}_d(M)$ such that $Q_N=Q_C$. 
	\end{lemma}
	
	\begin{proof}
		For each choice of $i$ and $j$, we define  $c_{ij}=\max\{u_i-u_j : u\in Q_N\}$. The matrix $C=(c_{ij})$ lives in $\RR_0^{d\times d}$ and has the property that $Q_N=Q_C$. Moreover, since $Q_N$ is contained in $Q_M$, we have  $C\leq M$. The fact  that $C\mintimes C=C$ follows from the definition of the $c_{ij}$'s and \eqref{eq:polytrope1}. In particular, $C$ belongs to the truncated polytrope region $\mathcal{P}_d(M)$.
	\end{proof}
	
	On the algebraic side,
	$\mathcal{P}_d(M)$ parametrizes all $\Ocal_K$-orders 
	$\Lambda_N$ that contain the given order $\Lambda_M$.
	Here $M$ and $N$ are assumed to be integer matrices.
	In particular, the integer points $u$ in $Q_M$ 
	correspond to maximal orders $\Lambda_{M(u)} =  \End_{\Ocal_K}(L_u) $
	that contain $\Lambda_M$; cf.~Proposition~\ref{prp:max order}.

	\begin{example}
		Let $M$ be the $d \times d$ matrix with entries $0$ on 
		the diagonal and $1$ off the diagonal.
		Thus $Q_M$ is the  pyrope \cite[\S 3]{joskul}.
		We consider two cases:
		the hexagon $(d=3)$ and Example \ref{ex:rhombic} $(d=4)$.
		The truncated polytrope region $\mathcal{P}_d(M)$
		classifies subpolytropes of~$Q_M$.
		
		\smallskip
		
		\noindent \underbar{$d=3$}:
		The $6$-dimensional polytope $\mathcal{P}_3(M)$ has the f-vector $(36,132,199,151,60,12)$. 
		Its $36$ vertices come in ten symmetry classes. We list the corresponding $3 \times 3$ matrices:		
		\[ \!\!\! \begin{small} \begin{matrix} 
				1, 6& \!\!\! [1, \!1, \! 1, \!1, \!1,\! 1]	&2, 6& \!\!\! [1, \! \frac{1}{2}, \!\frac{1}{2},\! 1,\! 1,\! \frac{1}{2}] 	&3, 8&  \!\!\! [0, \!-1\!, 0,\! -\!1, 1, 1] &
				3, 8& \!\!\! [1, 0, \!-\!1,\! -\!1, 0, 1] 	& 3, 8& \!\!\! [1, \!0, \!1,\! 1,\! 0,\! 1]  	 \\
				3, 6& \!\!\! [1, \! 1,\! 1,\! 1,\! 0,\! 0] &
				3, 6& \!\!\! [0,\! 1,\! 1,\! 1,\! 1,\! 0] 	&6, 7& \!\!\! [0, -1, 1, 0, 1, 1]  	&6, 7& \!\!\! [1, 1, 1, 1, 0, 1] &
				6, 6& \!\!\! [0, \! 0,\! 1,\! 1,\! 1,\! 0]
		\end{matrix} \end{small} \!\!\!\]
		These  polytropes are 
		shown in red in Figure~\ref{fig:eightpolytropes}. Our
		classification into $S_3$-orbits is finer than that  from
		symmetries of the hexagon $Q_M$, which leads to only eight orbits. 
		For us, this classification is more natural 
		because it reflects algebraic properties of orders.
		It distinguishes 
		min-plus vertices from max-plus vertices of $Q_M$.
		The polytope $\mathcal{P}_3(M)$ has $41$ integer points, so there are 
		$41$ orders containing~$\Lambda_M$. In addition to $34$ integer vertices, there are seven interior integer points, namely $[0,0,0,0,0,0]$ and  six like $[0,0,0,0,1,1]$,
		not seen in Figure~\ref{fig:eightpolytropes}.
		
		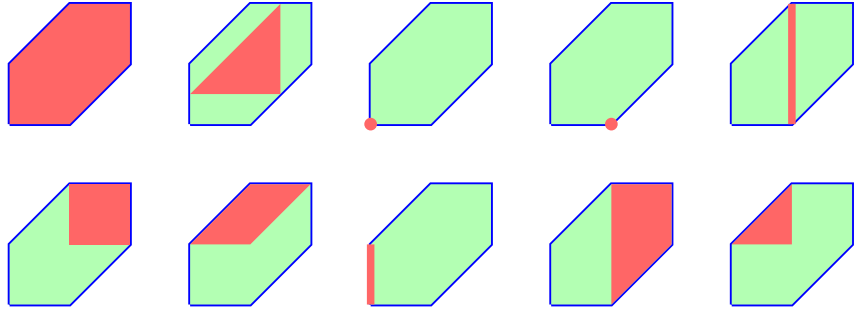
\begin{figure}[h]
			\begin{center}
				\vspace{-0.1cm}				
				\begin{tikzpicture}[scale=0.8]
					\draw[thick, blue, line width = 0.5mm]  (-1,-1) -- (0,-1) -- (1,0) -- (1,1) -- (0,1) -- (-1, 0) -- (-1 , -1);
					\fill[red!60]  (-1,-1) -- (0,-1) -- (1,0) -- (1,1) -- (0,1) -- (-1, 0) -- (-1 , -1);
					
					\draw[thick, blue, line width=0.5mm]  (2,-1) -- (3,-1) -- (4,0) -- (4,1) -- (3,1) -- (2, 0) -- (2 , -1);
					\fill[green!30]  (2,-1) -- (3,-1) -- (4,0) -- (4,1) -- (3,1) -- (2, 0) -- (2 , -1);
					\fill[red!60]  (3.5,1) -- (2, -0.5) -- (3.5, -0.5) -- (3.5, 1);
					
					\draw[thick, blue, line width=0.5mm]  (5,-1) -- (6,-1) -- (7,0) -- (7,1) -- (6,1) -- (5, 0) -- (5 , -1);
					\fill[green!30]  (5,-1) -- (6,-1) -- (7,0) -- (7,1) -- (6,1) -- (5, 0) -- (5 , -1);
					\fill[red!60] (5,-1) circle (3pt) node[below left= 1.5mm]{};				
					
					\draw[thick, blue, line width = 0.5mm]  (8,-1) -- (9,-1) -- (10,0) -- (10,1) -- (9,1) -- (8, 0) -- (8 , -1);
					\fill[green!30]   (8,-1) -- (9,-1) -- (10,0) -- (10,1) -- (9,1) -- (8, 0) -- (8 , -1);
					\fill[red!60] (9,-1) circle (3pt) node[below left= 1.5mm]{};
					
					\draw[thick, blue, line width = 0.5mm]  (11,-1) -- (12,-1) -- (13,0) -- (13,1) -- (12,1) -- (11, 0) -- (11 , -1);
					\fill[green!30]   (11,-1) -- (12,-1) -- (13,0) -- (13,1) -- (12,1) -- (11, 0) -- (11 , -1);
					\draw[thick, red!60, line width = 1mm]  (12,1) -- (12,-1);


					\draw[thick, blue, line width=0.5mm]  (-1,-4) -- (0,-4) -- (1,-3) -- (1,-2) -- (0,-2) -- (-1, -3) -- (-1 , -4);
					\fill[green!30] (-1,-4) -- (0,-4) -- (1,-3) -- (1,-2) -- (0,-2) -- (-1, -3) -- (-1 , -4);
					\fill[red!60]  (0,-3) -- (1,-3) -- (1,-2) -- (0,-2)--(0,-3);

					\draw[thick, blue, line width = 0.5mm]  (2,-4) -- (3,-4) -- (4,-3) -- (4,-2) -- (3,-2) -- (2, -3) -- (2 , -4);
					\fill[green!30]  (2,-4) -- (3,-4) -- (4,-3) -- (4,-2) -- (3,-2) -- (2, -3) -- (2 , -4);
					\fill[red!60]     (2,-3) --  (3,-2) -- (4,-2) -- (3,-3)-- (2,-3) ;
					
					\draw[thick, blue, line width = 0.5mm]   (5,-4) -- (6,-4) -- (7,-3) -- (7,-2) -- (6,-2) -- (5, -3) -- (5 , -4);
					\fill[green!30]    (5,-4) -- (6,-4) -- (7,-3) -- (7,-2) -- (6,-2) -- (5, -3) -- (5 , -4);
					\draw[thick, red!60, line width = 1mm]  (5,-3) -- (5,-4);
					
					\draw[thick, blue, line width=0.5mm]  (8,-4) -- (9,-4) -- (10,-3) -- (10,-2) -- (9,-2) -- (8, -3) -- (8 , -4);
					\fill[green!30]  (8,-4) -- (9,-4) -- (10,-3) -- (10,-2) -- (9,-2) -- (8, -3) -- (8 , -4);
					\fill[red!60]  (9,-4) -- (9,-2) -- (10,-2) -- (10,-3)--(9,-4);
					
					\draw[thick, blue, line width = 0.5mm]  (11,-4) -- (12,-4) -- (13,-3) -- (13,-2) -- (12,-2) -- (11, -3) -- (11 , -4);
					\fill[green!30]    (11,-4) -- (12,-4) -- (13,-3) -- (13,-2) -- (12,-2) -- (11, -3) -- (11 , -4);
					\fill[red!60]     (11,-3) --  (12,-3) -- (12,-2) -- (11,-3);
				\end{tikzpicture}
				\vspace{-0.25cm}
			\end{center}        
			\caption{The regular hexagon has $36$ extreme subpolytropes in ten symmetry classes.
				\label{fig:eightpolytropes}}
		\end{figure}     
		
		\noindent \underbar{$d=4$}:	
		The truncated polytrope region $\mathcal{P}_4(M)$
		for (\ref{eq:J4}) is $12$-dimensional.
		Its f-vector~is
		\[
		(961, 17426, 103780, 304328, 517293, 549723, 377520, 168720, 48417, 8620, 894, 48
		).
		\]
		The $961$ vertices come in $65$ orbits under the $S_4$-action. Among the simple vertices
		we find:
		\[ \begin{small} \begin{matrix}
				1,   12 & \!\!\!  [1, 1, 1, 1, 1, 1, 1, 1, 1, 1, 1, 1] \quad & \qquad  &
				8,   12& \!\!\! [1, 1, 1, 1, \frac{1}{2}, 1, 1, 1, \frac{1}{2}, 1, \frac{1}{2}, 1] 
				\\
				4,   27& \!\!\!  [1, 1, 1, -1, 0, 0, -1, 0, 0, -1, 0, 0] \quad &   \qquad &
				4,   27& \!\!\!  [-1, -1, -1, 1, 0, 0, 1, 0, 0, 1, 0, 0]\\
		\end{matrix} \end{small}  \]		
		The list  of all vertices, and much more, is made available at 
		our {\tt mathrepo} site.
		Such data sets can be useful for comprehensive computational studies of $\Ocal_K$-orders in $K^{d \times d}$.
	\end{example}
	
	\section{Ideals} \label{sec5}
	
	To better understand the order $\Lambda_M$ for $M \in \mathcal{P}_d$, we study its (fractional) ideals. By an \emph{ideal} of $\Lambda_M$ we mean an additive subgroup $I$ of $\Lambda_M$ such that $\Lambda_M I \subseteq I$ and $I \Lambda_M \subseteq I$.
	A {\em fractional ideal} of  $\Lambda_M$ is a (two sided) $\Lambda_M$-submodule $J$ of $K^{d \times d}$ such that $\alpha J \subset \Lambda_M$ for some $\alpha \in K^*$.
	
	\begin{example} \label{ex:frac-ideals}
		Fix $X \in K^{d \times d}$ and consider  the two-sided $\Lambda_M$-module  
		$\langle X \rangle \,=\, \Lambda_M X \Lambda_M \,= \,\big\{ A X B\, \colon 
		A,B \in \Lambda_M \big\}.$
		This is an ideal when $X \in \Lambda_M$.
		If $X \not\in \Lambda_M$ then $\alpha X \in \Lambda_M$ for some  $\alpha \in K^*$.
		Hence, $\langle X \rangle$ is a fractional ideal.  These are the {\em principal (fractional) ideals} of $\Lambda_M$.
	\end{example}
	
	For all that follows, we assume that $M \in \mathcal{P}_d$ is an integer matrix in standard form.
	
	\begin{proposition} \label{prop:I_n}
		The nonzero fractional ideals of the order $\Lambda_M$ are the sets of the form
		\begin{equation} \label{eq:I_n}
			I_N \,\, = \,\, \bigl\{ \,X \in K^{d \times d} \colon \val(X) \geq N  \,\bigr\}, 
		\end{equation}
		where $N = (n_{ij}) $ is any matrix in $ \,\ZZ^{d \times d }$ with
		$ N \mintimes M = 	M \mintimes N = N$. This is equivalent to
		\begin{equation}\label{frac_ideal_ineqs}
			n_{ik} \leq n_{ij}  + m_{jk}  \quad {\rm and} \quad
			n_{ik} \leq m_{ij}  + n_{jk} \qquad {\rm for} \quad
			1 \leq i,j,k \leq d.
		\end{equation}
	\end{proposition}
	
	\begin{proof}
		The result is due to Plesken who states it in (viii) from
		\cite[Remark~II.4]{plesken83}.
		The min-plus matrix identity $N \mintimes M = N$ is equivalent to
		$n_{ik} \leq n_{ij}  + m_{jk} $ because $m_{jj} = 0 $.
	\end{proof}
	
	\begin{remark}
		If $N$ has zeros on its diagonal 
		and satisfies (\ref{eq:polytrope1})
		then $I_N = \Lambda_N$ is an order, as before.
		However, among
		all lattices in $K^{d \times d}$, ideals are more general than orders.
		In particular, we generally have $n_{ii} \not= 0$ for the matrices $N$
		in (\ref{eq:I_n}).
		A fractional ideal $I_N$  is an ideal in $\Lambda_M$
		if and only if  $N \geq M$. If this holds then the polytrope $Q_N$ is contained in~$Q_M$.
	\end{remark}

	\begin{example}
		The Jacobson radical of the order $ \Lambda_M$ is 
		the ideal ${\rm Jac}(\Lambda_M) = I_{M+{\rm Id}_d}$.
		Here ${\rm Id}_d$ is the identity matrix. The quotient of $\Lambda_M$ 
		by its Jacobson radical is the product of residue fields
		$\Lambda_M/{\rm Jac}(\Lambda_M)\cong 
		(\Ocal_K/\langle p \rangle)^d$. See (i) in \cite[Remark~II.4]{plesken83} for more details.
	\end{example}

	Let $\mathcal{Q}_M$ denote the set of matrices $N$ in $\RR^{d \times d}$ that satisfy the inequalities in (\ref{frac_ideal_ineqs}). These inequalities are bounds on
	differences of matrix entries in $N$.
	We can thus regard $\mathcal{Q}_M$ as a polytrope in $\RR^{d \times d} / \RR {\bf 1}$, where ${\bf 1} = \sum_{i,j=1}^d E_{ij}$. The matrices $N$ parameterizing  the fractional ideals
	$I_N$ of $\Lambda_M$ (up to scaling) are the integer points of $\mathcal{Q}_M$. One checks directly that $\mathcal{Q}_M$ is closed under both addition and multiplication of matrices in the min-plus algebra. Its product $\mintimes$ represents the multiplication of fractional ideals as the following proposition~shows.
	
	\begin{proposition}
		If $M \in \mathcal{P}_d$ is in standard form and
		$N,N' \in \mathcal{Q}_M$ then $\,I_{N} I_{N'}\, =\, I_{N \mintimes N'}$.
	\end{proposition}
	
	\begin{proof} Let $X \in I_N,  Y\in I_{N'}$. The inequalities $\val(X) \geq N,\,\val(Y) \geq N'$ imply
		$\val(XY) \geq \val(X) \mintimes \val(Y) \geq N \mintimes N'$ and so
		$X Y \in I_{N \mintimes N'}$. This gives the inclusion $I_N I_{N'} \subseteq I_{N \mintimes N'}$.
		Let $u_{ij} = \min\limits_{1 \leq k \leq d} (n_{ik} + n'_{kj})$ be the $(i,j)$ entry of $N \mintimes N'$.
		For the inclusion $ I_{N \mintimes N'} \subseteq I_N I_{N'} $,
		it suffices to show that $p^{u_{ij}} E_{ij} $ is in $ I_N I_{N'}$ for all $i,j$.
		Fix $i,j$ and let $k$ satisfy $u_{ij} = n_{ik} + n'_{kj}$. 
		The matrices $p^{n_{ik}} E_{ik}$ and $p^{n_{kj}'} E_{kj}$ are in $I_N$ and $I_{N'}$.
		Their product $p^{u_{ij}} E_{ij}$ is in $I_{N} I_{N'}$.
	\end{proof}	
	
	We call $\mathcal{Q}_M$ the {\em ideal class polytrope} of $M$.
	The min-plus semigroup $(\mathcal{Q}_M, \mintimes)$ plays the role of the ideal class group in number theory. Its neutral element is 
	the given matrix $M$.

	\begin{example} \label{ex:octahedron}
		Fix $M =  \begin{bmatrix} 0  &  1\\ 1 & 0 \end{bmatrix} \in \mathcal{P}_2$. 
		The polytrope $\mathcal{Q}_M$ is the octahedron with vertices
		\[ 
		\begin{small}
			\begin{bmatrix}	
				0 & 1 \\
				1 & 2 
			\end{bmatrix} ,\,\,
			\begin{bmatrix}	
				2 & 1 \\
				1 & 0 
			\end{bmatrix} ,\,\,\,\,\begin{bmatrix}	
				1 & 2 \\
				0 & 1 
			\end{bmatrix} ,\,\,\begin{bmatrix}	
				1 & 0 \\
				2 & 1 
			\end{bmatrix} ,\,\,\,\,\begin{bmatrix}	
				0 & 1 \\
				1 & 0 
			\end{bmatrix} ,\,\,\begin{bmatrix}	
				1 & 0 \\
				0 & 1 
			\end{bmatrix} 
		\end{small}
		\quad \,\in \,\, \ZZ^{2 \times 2}/\ZZ {\bf 1}.
		\]
		This octahedron contains  $19$ integer points $N$. These are in bijection with the 
		equivalence classes of fractional ideals $I_N$ in the order $\Lambda_M$.
		The midpoint of $\mathcal{Q}_M$ corresponds to the Jacobson radical  $I_{M+{\rm Id}_2}$.
		The remaining $12$ integer points are the midpoints of the edges.
	\end{example}
	
	One may ask whether the ideal class semigroup $(\mathcal{Q}_M, \mintimes)$ is actually a group. To address this question, we define the {\em pseudo-inverse} of a fractional ideal $I$
	in the order $\Lambda_M$  as follows:
	\[ (\Lambda_M: I) \,\,= \,\,\{ \,X\in K^{d \times d} \,\colon XI \subseteq \Lambda_M \,
	\textrm{ and } \, IX \subseteq \Lambda_M\}. \]
	
	\begin{lemma} \label{lem:inverse}
		The pseudo-inverse of a fractional ideal in $\Lambda_M$ is a fractional ideal in~$\Lambda_M$.
	\end{lemma}
	
	\begin{proof}
		Let $A \in \Lambda_M$ and $X \in (\Lambda_M : I)$, so that
		$ XI, IX \subseteq \Lambda_M$.
		Since $I$ is a fractional ideal,
		we have $AI \subseteq I$ and $IA \subseteq I$. 
		From these inclusions
		we deduce that $XA I , IXA, AXI , IAX$ are all subsets of $
		\Lambda_M$.  This implies $XA, AX \in (\Lambda_M : I)$.
		Hence $(\Lambda_M:I)$ is a fractional~ideal.	\end{proof}
	
	\begin{proposition}\label{prop:inverse}
		Let $M \in \mathcal{P}_d$ in standard form and $N \in \mathcal{Q}_M$.
		Then $(\Lambda_M: I_N) = I_{N'}$ where
		\begin{equation}\label{eq:inv} \quad
			n'_{ij}\,\, =\,\, \max\limits_{ 1 \leq \ell \leq d } \bigl(\, 
			\max (m_{\ell j} - n_{\ell i}, m_{i\ell} - n_{j\ell}) \, \bigr)
			\qquad {\rm for} \,\,\, 1 \leq i , j \leq d.
		\end{equation}
	\end{proposition}
	
	\begin{proof} By Proposition \ref{prop:I_n} and Lemma \ref{lem:inverse},
		there exists $N' \in \mathcal{Q}_M$ such that $I_{N'} = (\Lambda_M \colon I_N)$.
		Then $I_{N'} I_{N} \subseteq \Lambda_M$ and $I_N I_{N'} \subseteq \Lambda_M$,
		and $I_{N'}$ is the largest fractional ideal with this property.
		These two conditions are equivalent to  $\,p^{n'_{ij}} E_{ij} I_N \subseteq \Lambda_M\,$ and
		$\,p^{n'_{ij}} I_N E_{ij}  \subseteq \Lambda_M $ for all $i,j$.
		The first condition holds if and only if $n'_{ij} + n_{j\ell} \geq m_{i\ell}$ for all $\ell$.
		The second condition holds if and only if  $n_{\ell i} + n'_{ij} \geq m_{\ell j}$ for all $\ell$.
		The smallest solution $N' = (n_{ij}')$ is given by (\ref{eq:inv}).
	\end{proof}
	
	Passing from ideals to their matrices, we also call $N'$ the
	{\em pseudo-inverse} of~$N$ in $\mathcal{Q}_M$.
	
	\begin{example} \label{ex:nineteen} Let $d=2$ and $M$ as in Example~\ref{ex:octahedron}.
		The $19$ ideal classes $N$ in $\mathcal{Q}_M$ have only
		three distinct pseudo-inverses:
		$\,N' \in  \bigl\{ 
		\begin{tiny} \begin{bmatrix} 
				0 \! & \! 0 \\ 
				0 \! & \! 0  \end{bmatrix} \end{tiny},\,
		\begin{tiny} \begin{bmatrix} 
				0 \! & \! 1 \\ 
				1 \! & \! 0  \end{bmatrix} \end{tiny},\,
		\begin{tiny} \begin{bmatrix} 
				1 \! & \! 0 \\ 
				0 \! & \! 1  \end{bmatrix} \end{tiny} \bigr\}$.
		For most ideal classes $N$, we have
		$N \mintimes N' \neq M$ and $ N' \mintimes N \neq M$.
		This means that most $N$ do not have an inverse in $(\mathcal{Q}_M, \mintimes)$.
		In particular, the ideal class polytrope  $\mathcal{Q}_M$ is a semigroup but not a group. 
	\end{example}
	
	The semigroup $\mathcal{Q}_M$ has the neutral element $M$ and
	each ideal class $N \in \mathcal{Q}_M$ has 
	a pseudo-inverse $N'$ given by the formula (\ref{eq:inv}).
	With this data, we define the
	\emph{ideal class group}  
	\[
	\mathcal{G}_M \,\,= \,\, \bigl\{ \, N \in \mathcal{Q}_M \,:\, 
	N \mintimes N' \,=\, N' \mintimes N \, =\, M \bigr\}.
	\]
	This is the maximal subgroup of the semigroup $\mathcal{Q}_M$. 
	It would be interesting
	to understand how $M$ determines the structure of $\mathcal{G}_M$.
	Note that
	$\,\mathcal{G}_M = \bigl\{ 
	\begin{tiny} \begin{bmatrix} 
			0 \! & \! 1 \\ 
			1 \! & \! 0  \end{bmatrix} \end{tiny},\,
	\begin{tiny} \begin{bmatrix} 
			1 \! & \! 0 \\ 
			0 \! & \! 1  \end{bmatrix} \end{tiny} \bigr\}\,$
	in Example \ref{ex:nineteen}. 
	
	\begin{example} \label{ex:ICG}
		Here are three examples of ideal class groups of graduated orders:
		$$  \begin{matrix}
			M_2 = \begin{tiny} \begin{bmatrix} 0  & 1 \\ 1 & 0 \end{bmatrix} \end{tiny} & &
			M_3 = \begin{tiny} \begin{bmatrix} 0 & 1 & 1 \\ 1 & 0 & 1 \\ 1 & 1 & 0 \end{bmatrix} \end{tiny}  & &
			M_4 = \begin{tiny} \begin{bmatrix} 0 & 1 & 1 & 1 \\ 1 & 0 & 1 & 1 \\ 1 & 1 & 0 & 1 \\ 1 & 1 & 1 & 0 \end{bmatrix} \end{tiny}& &\\
			\mathcal{G}_{M_2} \cong  \ZZ/2\ZZ \phantom{ooo} &  &
			\mathcal{G}_{M_3} \cong  \ZZ / 6\ZZ \phantom{ooo} &  &
			\mathcal{G}_{M_4} \cong  S_4 \phantom{ooo} & &
		\end{matrix}
		$$
		The isomorphism types of these groups were computed using {\tt GAP}; the code is at
		our {\tt mathrepo} site.
We do not know
		how this list continues 
		for pyropes
		\cite[\S 3]{joskul}
		in higher dimensions.
	\end{example}		
	
	We end this section with a conjecture about the geometry of
	$\mathcal{G}_M$ inside $\mathcal{Q}_M$.

	\begin{conjecture}
		For any integer matrix $M $ in the polytrope region $ \mathcal{P}_d$, the elements in the 
		ideal class polytrope $\mathcal{G}_M$ are among the classical vertices of  
		the ideal class polytrope $\mathcal{Q}_M$.
	\end{conjecture}

	\section{Towards the building}		
	\label{sec6}
	
	Affine buildings \cite{brown, zhang} provide a natural setting for orders and min-max convexity.
	The objects we discussed in this paper so far are associated to one
	apartment in this building, namely, that corresponding to the diagonal lattices.
	The aim of this section is to present
	this perspective and to lay the foundation for a general theory that
	goes beyond one apartment.
	
	\begin{definition}
		The \emph{affine building} $\Bcal _d(K)$ is an infinite simplicial 
		complex. Its vertices  are the equivalence classes $[L]$ of lattices in $K^d$.
		A configuration $\{ [L_1], \ldots , [L_s] \} $ is a simplex 
		in $\Bcal _d(K)$ if and only if, up to some permutation, there exist representatives 
		$\tilde{L}_i \in [L_i]$ satisfying
		$\tilde{L}_1 \supset \tilde{L}_2 \supset \cdots \supset \tilde{L}_s \supset p \tilde{L}_1$.
		The maximal simplices $\{ [L_1] ,\ldots , [L_d] \} $ are called \emph{chambers}.
		The {\em standard chamber} $C_0 $ is given by the diagonal lattices
		$ L_i = L_{({\bf 1}_{i-1}, {\bf 0}_{d-i+1})} = L_{(1,\dots,1,0,\dots,0)}$.
	\end{definition}
	
	Given a basis $\{b_1, \dots, b_d\}$ of $K^d$, the \emph{apartment} defined by this basis is the 
	set of classes $[L]$ of all lattices $L= \bigoplus _{i=1}^d p^{u_i} \Ocal _K b_i $ 
	where $u_1,\ldots , u_d $ range over $ \ZZ $. Hence the apartment is
	\[
	\big\{ \left[  p^{u_1} \Ocal_K b_1 \oplus \dots \oplus p^{u_d} \Ocal_K b_d \right] \colon u_1, \dots, u_d \in \ZZ \,\big\} \,\,=\,\, \bigl\{\, [gL_u] \,\colon u \in \ZZ^d \,\bigr\},
	\]
	where $g \in \GL_d(K)$ is the matrix with columns $b_1, \dots, b_d$.
	The {\em standard apartment} is 
	the one associated with the standard basis $(e_1,\dots, e_d)$ of $K^d$. 
	The vertices of the standard apartment are the diagonal lattice classes $[L_u]$ for $u \in \ZZ^{d}$.
	We identify this set of vertices with $\ZZ^n / \ZZ {\bf 1}$.
	
	The general linear group $\GL _d(K)$ acts on the building $\Bcal _d(K)$. This  action
	preserves the simplicial complex 
	structure. In fact, the action is transitive on lattice classes, on apartments and also on the chambers. 
	The stabilizer of the standard lattice $L_0$ is the 
	subgroup 
	\[\GL_d(\Ocal _K) \,\,= \,\, \{ \,g\in \Ocal _K^{d\times d} \,:\, \val(\det(g)) = 0 \,\}
	\,\, \subset \,\, \GL_d(K).
	\]
	Starting from the standard chamber $C_0$, there exist reflections 
	$s_0,s_1,\ldots , s_{d-1}$ in $\GL_d(K)$ that map $C_0$ to the $d$ adjacent 
	chambers in the standard apartment.  For $i \geq 1$, define 
	$s_i$ by
	\[
	s_i (e_i) = e_{i+1} \text{,} \quad s_i(e_{i+1}) = e_i \, \text{ and } \, s_{i}(e_j) = e_j \text{ when } j \neq i, i+1.
	\]
	The map $s_0$ is defined by $s_0(e_i ) = e_i$ for $i=2,\ldots , d-1$ 
	and $s_0(e_d) = p e_1$, $s_0(e_1) = p^{-1} e_d$. 
	The reflections $s_0,\ldots , s_{d-1}$ are Coxeter generators for the {\em affine Weyl group}
	$\,W = \langle s_0 ,\ldots , s_{d-1 }\rangle$.
	The group $W$ acts regularly on the chambers $C$ in the standard apartment 
	\cite[\S~1.5, Thm.~2]{bourbaki}: for every $C$ there is a unique
	$w\in W$ such that $C=wC_0$. The elements of $W$ are the matrices $h_{\sigma} g_u$ where $h_{\sigma} = (1_{i = \sigma(j)})_{i,j}$ for $\sigma \in S_d$,  and
	$u \in \ZZ^{d}$ with $ u_1+\cdots+u_d = 0$. Thus $W$ is the semi-direct product of $S_d$ and the group of diagonal matrices $g_u$ whose exponents sum to $0$.
	
	Our primary object of interest is the Plesken-Zassenhaus order ${\rm PZ}(\Gamma)$
	of a finite configuration $\Gamma$ in the affine building $\Bcal_d(K)$.
	This is the intersection (\ref{eq:PZ1}) of endomorphism~rings.
	In this paper we studied the case when $\Gamma$ lies in one apartment.
	In Theorem \ref{thm:PZtropical} we showed that ${\rm PZ}(\Gamma) = \Lambda_M$
	where $M$ is the matrix in $\mathcal{P}_d$ that encodes 
	the min-max convex hull of $\Gamma$.
	This was used in Sections \ref{sec4} and \ref{sec5} to
	elucidate combinatorial and algebraic structures in ${\rm PZ}(\Gamma)$.
	A subsequent project will extend our results to arbitrary configurations $\Gamma$
	in $\Bcal_d(K)$. 
	
	We conclude this article with configurations
	given by two chambers $C,C'$ in $\mathcal{B}_d(K)$.
	We are interested in the their order $\PZ(C \cup C')$.
	A fundamental fact about buildings
	states that any two chambers $C,C'$ lie
	in a common apartment, cf.~\cite{bourbaki, brown}.  Also, since the affine Weyl group $W$ acts regularly on the chambers of the standard apartment, we can then reduce to the case where the two chambers in question are $C_0$ and $wC_0$ for some $w = h_{\sigma} g_u \in W $.

	\begin{example}\label{ex:stdChamber}
		The standard chamber $C_0$ is encoded by 
		$M_0 = \sum_{1 \leq i < j \leq d} E_{ij}$. The polytrope $Q_{M_0}$ is a simplex.
		The order ${\rm PZ}(C_0) = \Lambda_{M_0}$ consists of all
		$X \in \Ocal_K^{d \times d}$ with $x_{ij} \in \langle p \rangle $ for~$i < j$.
	\end{example}

	Let $D_u = \val(g_u)$ denote the tropical diagonal matrix with  $u_1, \dots, u_d$ 
	on the diagonal and $+\infty$ elsewhere. We also write $P_\sigma \coloneqq \val(h_\sigma)$
	for the tropical permutation matrix given by $\sigma$.
	
	\begin{proposition}
		We have $ {\rm PZ}(C_0 \cup h_{\sigma}g_u C_0) = \Lambda_{M^{\sigma,u}}$
		where the matrix $M^{\sigma,u}$ is given by
		\[
		M^{\sigma,u} \,\,=\,\, M_0\, \,\maxplus\,
		\left(  P_\sigma \,\mintimes\, D_{u} \,\mintimes \,
		M_0 \,\mintimes \, D_{-u} \,\mintimes \, P_{\sigma^{-1}}  \right).
		\]
	\end{proposition}
	
	\begin{proof}
		We have $\PZ(C_0 \cup h_{\sigma}g_u C_0) = \PZ(C_0) \cap \PZ(h_{\sigma}g_u C_0)$. Recall that $\PZ(C_0) = \Lambda_{M_0}$ from \cref{ex:stdChamber}.
		Suppose that $M \in \ZZ_{0}^{d \times d}$ satisfies $\PZ(h_{\sigma}g_u C_0) = \Lambda_{M}$. 
		By Theorem \ref{thm:PZtropical},	 the order  $\Lambda_{M_0 \maxplus M}$ is equal to
		$\PZ(C_0 \cup h_{\sigma}g_u C_0)$. To determine $M$, notice that $\PZ(wC_0) = h_\sigma g_u \PZ(C_0) g_{-u} h_{\sigma^{-1}} $. This implies
		the stated formula $\,M = P_\sigma \mintimes D_{u} \mintimes  M_0 \mintimes D_{-u} \mintimes P_{\sigma^{-1}} $.
	\end{proof}
	
	We may ask for invariants of the orders $\PZ(C_0 \cup w C_0)$ in terms of $w \in W$.
	Clearly, not all polytropes in an apartment arise as the min-max convex hull of two chambers. 
	Which graduated orders are of the form $\PZ(C_0 \cup wC_0) $? 
	Which other elements $w ' $ in the affine Weyl group $W$ give rise to the same
	Plesken-Zassenhaus order $\PZ(C_0 \cup wC_0) $ up to isomorphism?
	
	\bigskip \bigskip
	
	{\bf Acknowledgements}.
	This project was supported by the Deutsche Forschungsgemeinschaft 
	(DFG, German Research Foundation) -- Project-ID 286237555 --TRR 195.
	We thank Michael Joswig for his comments on an early version of this paper.

	\bigskip
	\bigskip
	
	\noindent
	\footnotesize 
	{\bf Authors' addresses:}
	
	\smallskip
	
	\noindent Yassine El Maazouz,
	UC Berkeley,
	\hfill {\tt yassine.el-maazouz@berkeley.edu}
	
	\noindent Marvin Hahn,  MPI-MiS Leipzig
	\hfill {\tt  marvin.hahn@mis.mpg.de}
	
	\noindent Gabriele Nebe, RWTH Aachen
	\hfill {\tt gabriele.nebe@rwth-aachen.de}
	
	\noindent Mima Stanojkovski,  RWTH Aachen and
	MPI-MiS Leipzig
	\hfill {\tt mima.stanojkovski@mis.mpg.de}
	
	\noindent Bernd Sturmfels,
	MPI-MiS Leipzig and UC Berkeley 
	\hfill {\tt bernd@mis.mpg.de}
	

\bigskip
	
	\begin{thebibliography}{alpha} 
		\begin{small} 
			\setlength{\itemsep}{-0.4mm}
			
						\bibitem{brown} Peter Abramenko and Kenneth Brown: {\em Buildings: Theory and Applications}, Springer-Verlag, New York, 2008.
						
			\bibitem{avis} David Avis: {\em On the extreme rays of the metric cone},	Canadian J.~Math.~{\bf 32} (1980) 126--144.
			
			\bibitem{bourbaki} Nicolas Bourbaki:
			{\em Lie groups and {L}ie algebras, {C}hapters 4--6}, Elements of Mathematics,
			Springer-Verlag, Berlin, 2002.
			
			\bibitem{CGQ04}	Guy Cohen and St\'{e}phane Gaubert and Jean-Pierre Quadrat: {\em Duality and separation theorems in idempotent semimodules}, 
			Linear Algebra Appl. {\bf 379} (2004) 395--422.
			
			\bibitem{deza} Antoine Deza: {\em Metric polytopes and metric cones}, \url{www.cas.mcmaster.ca/~deza/metric.html}.

\bibitem{DDP}
Michel~Deza, Mathieu~Dutour and Elena~Panteleeva:
{\em Small cones of oriented semi-metrics},
 Amer. J. Math. Management Sci. {\bf 22} (2002) 199--225. 
 
 \bibitem{DKKP}
 Mikhailo Dokuchaev,  Vladimir Kirichenko, Ganna Kudryavtseva and Makar Plakhotnyk:
 {\em The max-plus algebra of exponent matrices of tiled orders}, J. Algebra {\bf 490} (2017) 1--20.
 
			
			\bibitem{polymake} Ewgenij Gawrilow and Michael Joswig:
			{\em Polymake: a framework for analyzing convex polytopes},
			Polytopes -- combinatorics and computation (Oberwolfach, 1997), 43--73, DMV Seminar {\bf 29}, Birkh\"auser, Basel, 2000.
			
			\bibitem{tiled} Vasanti Jategaonkar:
			{\em Global dimension of tiled orders over a discrete valuation ring},
			Trans. Am. Math. Soc. {\bf 196} (1974) 313--330. 
			
			\bibitem{joswig} Michael Joswig: {\em Essentials of Tropical Combinatorics},
			Graduate Studies in Mathematics, American Mathematical Society, 2022.
			
			\bibitem{joskul} Michael Joswig and Katja Kulas:
			{\em Tropical and ordinary convexity combined},
			Adv. Geom. {\bf 10} (2010) 333--352.
			
			\bibitem{maclagan} Diane Maclagan and Bernd Sturmfels:
			{\em Introduction to Tropical Geometry}, Graduate Studies in 
			Mathematics, Vol 161, American Mathematical Society, 2015. 
			
			\bibitem{plesken83} Wilhelm Plesken: {\em Group Rings of Finite Groups
				over $p$-adic Integers}, Lecture Notes in Mathematics, {\bf 1026}, Springer-Verlag, Berlin, 1983.
			
			\bibitem{shemanske} Thomas Shemanske: {\em Split orders and convex polytopes in buildings}, J. Number Theory {\bf 130} (2010) 101--115.
			
			\bibitem{tran} Ngoc Tran: {\em Enumerating polytropes}, J. Combin. Theory, Ser.~A {\bf 151} (2017)~1--22.
			
			\bibitem{monomial} 
			Tse-Chung Yang and  Chia-Fu Yu:
			{\em Monomial, Gorenstein and Bass orders},
			J. Pure Appl. Algebra {\bf 219} (2015) 767--778.
			
			\bibitem{zhang} Leon Zhang: {\em Computing min-convex hulls in the affine
				building of ${\rm SL}_d$}, Discrete Comput. Geom. {\bf 65} (2021) 1314--1336.
			
			
		\end{small}                                                        
	\end{thebibliography}
\end{document}